\documentclass[12pt,draftcls,onecolumn]{IEEEtran}
\usepackage{mathrsfs}
\usepackage{amsmath,amssymb}
\usepackage{graphicx}

\def\setzero{\setcounter{equation}{0}}

\def\Pr{{\noindent\bf Proof. }}

\begin{document}
\begin{center}
{\Large \bf Optimal control of uncertain stochastic systems with
Markovian switching  and its applications to portfolio
decisions}\footnote{Project supported by National Natural Science
Foundation of China (71171003) and Anhui Natural Science Foundation
of Universities
(KJ2012B019, KJ2013B023). \\
$^*$Corresponding author: wyfei@dhu.edu.cn. }

\vskip12pt

{\rm Weiyin Fei}$^*$

\small{ (School of Mathematics and Physics, Anhui Polytechnic
University, Wuhu 241000, Anhui, P.R. China)}
\end{center}

\vskip12pt \noindent{\bf Abstract:} This paper first describes
a class of uncertain stochastic control systems with  Markovian switching, and derives
 an It\^o-Liu formula for Markov-modulated processes. And we characterize an optimal control law,
which satisfies the generalized
Hamilton-Jacobi-Bellman (HJB) equation with Markovian switching.
Then, by using the generalized HJB equation, we deduce the optimal
consumption and portfolio policies under uncertain stochastic financial markets with Markovian switching.
Finally, for constant relative risk-aversion (CRRA) felicity functions, we explicitly obtain the optimal
consumption and portfolio policies. Moreover, we also make an economic
analysis through numerical examples.

\vskip12pt


\vskip12pt

%{\noindent\bf Mathematics Subject Classification (2000)}: 60H10; 60J27; 90C39;  91B28; 93E20

\vskip12pt

\section{Introduction} \setzero

\vskip12pt

Uncertainty traditionally contains randomness, fuzziness and uncertainty resulting from their actions together. In many cases, fuzziness and randomness
 simultaneously appear in a system. To describe this phenomenon, the concept of fuzzy random variables was introduced by Kwakernaak \cite{Kw} as a random variable taking fuzzy variable values. According to different requirements of measurability,
different studies for fuzzy random variables were made by Puri and Ralescu \cite{PR},  Liu and Liu \cite{LL}. The fuzzy random differential equations were studied by Fei \cite{FeiINS2007, FeiKBS}, and the references therein.

For stochastic dynamic systems, stochastic differential equations  were studied by such
researchers as Karatzas and Shreve \cite{KS1}.
When a fuzzy dynamic system is considered, fuzzy differential equations were investigated by Kaleva \cite{Ka}, Fei \cite{FeiF, FeiI}, etc.

 In a real world,  some information and knowledge are usually represented by human linguistic expressions
like ``about 100km", ``approximately 80kg", ``warm", ``young". Perhaps they are often considered to be
subjective probability or fuzzy concepts. However, much research showed that those ``unknown constants" and ``unsharp concepts" behave
neither like randomness nor like fuzziness. In order to distinguish this
phenomenon from randomness and fuzziness, Liu \cite{Liu4th} names it ``uncertainty".
After further observations, we find that randomness and Liu's uncertainty often appear simultaneously in a system.
To describe this phenomenon, the concept of uncertain random variable
  was introduced by Liu \cite{Liu2013} as a random variable
taking ``uncertain variable" values.

Under the stochastic context with Markovian switching, the control systems have
been used to model many practical systems where they may experience
abrupt changes in their structure and parameters caused by such phenomena as component failures or repairs, changing subsystem
interconnections, and abrupt environmental disturbances in an
economic system, etc. The control systems combine a part of the
state that is driven by the canonical process and Brownian motion and another part of the state
that takes discrete values.  Mariton \cite{Mar} explained that
the hybrid systems had been emerging as a convenient mathematical
framework for the formulation of various design problems in
different fields like target tracking (evasive target tracking
problem), fault tolerant control and manufacturing processes.

One of the important classes of the hybrid systems is the stochastic
differential equations with Markovian switching. In operation,
the system will switch from one mode to another in a random way, and
the switching between the modes is governed by a Markovian chain.
The optimal regulator, controllability, observability, stability and
stabilization etc. need to be studied. For more information about the
hybrid systems, the reader can refer to Feng et al. \cite{FLS}, Huang
and Mao \cite{HM}, and the references therein.

Hamilton \cite{Ha} originally proposed the regime switching models of stock
returns which  proved to be a better representation of financial
reality than the usual models with deterministic coefficients.
Thereafter, regime switching models have been studied in different
contexts. Option pricing in financial markets with regime switching
has been studied, for instance, by Buffington and Elliott \cite{BE},
 Guo and Zhang \cite{GZ2004}. In financial markets with regime switching, the
maximization of expected utility from consumption and/or terminal
wealth has also been studied by Sotomayor and Cadenillas
\cite{SC}, and Zhang and Yin \cite{ZQY},
etc. Under a mean-variance criterion, Zhou and Yin \cite{ZY}
propose a continuous-time Markowitz's mean-variance portfolio
selection model with regime switching which obtains the efficient
portfolio that minimizes the risk of terminal wealth given a fixed
expected terminal wealth. In Fei \cite{FeiTAC, FeiS}, the optimal control of Markovian switching systems with applications
to optimal consumption and portfolio under inflation and Markovian switching is dicussed. The optimal portfolio selection with random
fuzzy returns is explored in Huang \cite{HX}.

%%%%%%%%%%%%%%%%%%%%%%%%%%%%%%%%%%%%%%%%%%%%%%

Since randomness and Liu's uncertainty simultaneously appear in the financial market, we begin to consider the uncertain stochastic
systems with Markovian switching. The uncertain stochastic system driven by the canonical process and Brownian motion will switch from one mode to another in a random way, and
the switching between the modes is governed by a Markovian chain. After introducing the concept of uncertain random processes which is different from the one of You \cite{You2007}, we will generalize It\^o-Liu formula, from which we deduce Hamilton-Jacobi-Bellman (HJB) equation or the optimality equation for the optimal control law. Next, by using the optimality equation for uncertain stochastic controls, we solve the
optimal portfolio and consumption decision problem of an investor, which is its financial application. We obtain the explicit expression for the optimal consumption and portfolio. Finally, we make an economic analysis of the particular case of
the constant relative risk-aversion (CRRA) utilities.

In fact, the initial contribution to consumption/investment problems
in continuous-time was done by Merton \cite{Merton1969, Merton1971}.
Other models of consumption-investment can be found in Fei
\cite{FeiINS,FeiSM}, Fei and Wu \cite{FW2002}, Karatzas and Shreve \cite{KS2}, etc. The traditional
financial modeling only deals with the continuously stochastic changing
dynamics of a financial market. However, the market behavior is
also affected by uncertain stochastic processes and a
finite-state continuous-time Markovian chain that represents the
uncertainty generated by the more steady market conditions. Based on the uncertain random calculus, the optimal consumption and portfolio strategies are obtained under the uncertain random context. Our model differs from the
above papers and also answers a different economic question: What is
the effect of both Markovian switching and uncertainty on the optimal
consumption and portfolio?

 The rest of the paper is organized as follows. Section 2
provides the general framework of the uncertain stochastic optimal control problem with
Markovian switching.  Generalizing an It\^o-Liu formula, we deduce an HJB
equation for later research in the subsequent section.  The policy
of optimal consumption and portfolio under  the uncertain random environment with Markovian switching is derived in Section 3. In Section 4, for
the case of the CRRA utility, the optimal policies are explicitly
given, and the numerical results of optimal policies together with
an economic analysis of the investor's behavior are made. Finally,
our concluding remarks are presented in Section 5.

\vskip12pt

\section{Generalized It\^o-Liu formula and equation of optimality} \setzero

For convenience, we give some useful concepts at first. Let $\Gamma$ be a nonempty
set and $\cal L$ a $\sigma$-algebra over $\Gamma$. Each element $\Lambda\in {\cal L}$ is called an event. A set function $\cal M$ defined on the $\sigma$-algebra over $\cal L$, which satisfies (i) Normality; (ii) Monotonicity; (iii) Self-Duality and (iv) Countable Subadditivity, is called an uncertain measure according to Liu \cite{Liu4th}. Let a probability space be $(\Omega,{\cal A}, {\mathbb P})$.

\vskip12pt

%%%%%%%%%%%%%%%%%%%%%%%
{\noindent \bf Definition 2.1.} (Liu \cite{Liu4th}) Let $\cal T$ be an index set and let $(\Gamma, {\cal L}, {\cal M})$ be an uncertainty space. An uncertain process $X(t,\gamma)$ is a measurable function from ${\cal T}\times(\Gamma,{\cal L}, {\cal M})$ to the set of real numbers, i.e., for each $t\in {\cal T}$ and any Borel set $B$
of real numbers, the set $\{X(t)\in B\}=\{\gamma\in \Gamma| X(t,\gamma)\in B\}$ is an event.

\vskip12pt

The concepts and properties of the canonical process and other uncertain processes refer to Chapter 9 Liu \cite{Liu4th}.  In what follows, we give the notation of uncertain random variables.

\vskip12pt

{\noindent \bf Definition 2.2.} (Liu and Li \cite{LL}) An uncertain random variable is a function $\xi$ from a probability space $(\Omega, {\cal A}, {\mathbb P})$
 to the set of uncertain variables such that ${\cal M}\{\xi(\omega)\in B\}$ is a measurable function of $\omega$ for any Borel set $B$ of $\mathbb R$.

\vskip12pt
We now give the concept of expected value of an uncertain random variable.

\vskip12pt

{\noindent \bf Definition 2.3.}  Let $\xi$ be an uncertain random variable. Then its expected value is defined by
 $$E[\xi]=E_P[E_U[\xi]]$$
 provided that the  right hand operations are well defined. Here, the operators $E_P$ and $E_U$ stand for probability expectation and uncertain expectation, respectively.

 \vskip12pt

 Obviously, if both $a$ and $b$ are constant, then $E[aC_t+bB_t]=0$, where $C_t$ and $B_t$ are a scalar canonical process and a Wiener process (Brownian motion), respectively.

 \vskip12pt
{\noindent \bf Definition 2.4.}  A hybrid process $X(t, \gamma, \omega)$ is called an uncertain stochastic process if for each $t\in {\cal T}$, $X(t)$ is an uncertain random variable. An uncertain stochastic process $X(t)$ is called continuous if the sample paths of $X(t)$ are all continuous
functions of $t$ for almost all $(\gamma, \omega)\in ({\Gamma, \Omega})$.

 \vskip12pt
{\noindent \bf Definition 2.5.} (It\^o-Liu integral) Let $X(t)=(Y(t),Z(t))$ be an uncertain stochastic process.  For any partition of
closed interval $[a, b]$ with $a = t_1 < t_2 < \cdots < t_{N+1} = b,$ the mesh is written as $\Delta=\max\limits_{1<i<N}|t_{i+1}-t_i|.$
Then the It\^o-Liu integral of $X(t)$ with respect to $(B_t, C_t)$ is defined as follows,
$$\int_a^bX(s) d(B_s, C_s)=\lim\limits_{\Delta\rightarrow0}\sum\limits_{i=1}^N(Y({t_i})(B_{t_{i+1}}-B_{t_i})+Z({t_i})(C_{t_{i+1}}-C_{t_i}))$$
provided that it exists in mean square and is an uncertain random variable, where $C_t$ and $B_t$ are one-dimensional canonical process and one-dimensional Wiener process, respectively.
 In this case, $X(t)$ is
called It\^o-Liu integrable. Specially, when $Y(t)\equiv0$, $X(t)$ is called Liu integrable.

\vskip12pt

Next, we deduce It\^o-Liu formula for the case of mult-dimensional uncertain stochastic processes.

\vskip12pt

{\noindent \bf Theorem 2.6.} {\sl Let $B=(B_t)_{0\leq t\leq T}=(B^1_t, \cdots,
B_t^m)^\top_{0\leq t\leq T}$ and $C=(C_t)_{0\leq t\leq T}=(C^1_t, \cdots, C_t^n)^\top_{0\leq t\leq T}$ be an $m$-dimensional
standard Wiener process and an $n$-dimensional canonical process, respectively.
Assume that uncertain stochastic processes $X_1(t), X_{2}(t), \cdots, X_{p}(t)$ are given by
$${d} X_{k}(t)=u_{k}(t){d}t+\sum\limits_{l=1}^mv_{kl}(t){d}B_{t}^l+\sum\limits_{l=1}^nw_{kl}(t){d}C_{t}^l,\quad k=1, \cdots, p, $$
where $u_{k}(t)$ are all absolute integrable uncertain stochastic processes, $v_{kl}(t)$ are all square integrable
uncertain stochastic processes and $w_{kl}(t)$ are all Liu integrable uncertain stochastic processes. For $k,l=1, \cdots, p$,  let $\frac{\partial G}{\partial t}(t, x_1, \cdots, x_p)$, $\frac{\partial G}{\partial x_k}(t, x_1, \cdots, x_p)$ and $\frac{\partial^2 G}{\partial x_kx_l}(t, x_1, \cdots, x_p)$ be continuously functions. Then we have
$$
\begin{array}{ll}
&d G(t, X_{1}(t), \cdots, X_{p}(t))\\
&=\frac{\partial G}{\partial t}(t, X_{1}(t), \cdots, X_{p}(t)){d }t+\sum\limits_{k=1}^p\frac{\partial G}{\partial x_k}(t, X_{1}(t), \cdots, X_{p}(t)){d}X_{k}(t)\\
&\quad+\frac{1}{2}\sum\limits_{k=1}^p\sum\limits_{l=1}^p\frac{\partial^2 G}{\partial x_k\partial x_l}(t, X_{1}(t), \cdots, X_{p}(t)){d}X_{k}(t){d}X_{l}(t),
\end{array}
$$
where ${d}B_{t}^k{d}B_{t}^l=\delta_{kl}{d}t, {d}B_{t}^k{d}t={d}C_{t}^\imath{d}C_{t}^\jmath={d}C_{t}^\imath{d}t={d}B_{t}^k{d}C_{t}^\imath=0,$ for
 $k,l=1,\cdots, m, \imath,\jmath=1, \cdots, n.$ Here
$$
\begin{array}{ll}
\delta_{kl}=\left\{\begin{array}{ll}
0,& {\rm if}\quad k\neq l\\
1, & \rm {otherwise}.
\end{array}
\right.
\end{array}
$$
}

\vskip12pt

{\Pr} Similar to the discussion  in You \cite{You2007}, we easily derive our Theorem 2.6 for uncertain stochastic processes. Thus, the proof is complete.
$ \Box$

\vskip12pt

{\noindent \bf Remark 2.7.} Similar to the discussions in You \cite{You2007}, we can obtain the counterparts of It\^o-Liu formulas for multi-dimensional uncertain stochastic processes. Here, the further discussions is omitted.

\vskip12pt

In reality, the uncertain stochastic systems may experience
abrupt changes in their structure and parameters, so we characterize this phenomenon by an uncertain stochastic context with Markovian switching.  To this end, we make some preliminaries.

%%%%%%%%%%%%%%%%%%%%%%

Now we consider
$m$-dimensional standard Brownian motion $(B_t)_{0\leq t\leq T}$ and $n$-dimensional canonical process $(C_t)_{0\leq t\leq T}.$
  We observe
a continuous-time, stationary, finite-state Markov chain
$\zeta=(\zeta(t))_{0\leq t\leq T}$, and denote by $\cal S$ the
state-space of this Markov chain: that is, for every $t\in[0,T],
\zeta(t)\in{\cal S}=\{1,2,\cdots, S\}$, where
$S\in\{1,2,3,\cdots\}$. Here, the stochastic process $\zeta(t, \omega)$ represents the regime (
mode) of the random environment at time $t$ and depends only on $\omega\in \Omega$. $S$ is the number of regimes.
Assume that the Markov chain has a strongly irreducible generator
$Q=(q_{ij})_{S\times S}$, where
$q_{ii}\stackrel{\triangle}{=}-\lambda_i<0$ and $\sum_{j\in{\cal
S}}q_{ij}=0$ for every regime $i\in{\cal S}$. Denote by ${\mathbb A}=\{{\cal A}_t, t\in[0,T]\}$ the $\mathbb P$-augmentation of filtration
$\{{\cal A}_t^{B,\zeta}\}$ generated by the stochastic processes
$B, \zeta$, where ${\cal A}_t^{B,
\zeta}\stackrel{\triangle}{=}\sigma\{B_s, \zeta(s),0\leq
s\leq t\}$ for every $t\in[0,T]$. We assume that the Markov chain
$\zeta(\cdot)$ is independent of the Brownian motion
$B$.

Let ${\mathbb U} $ be a separable metric space, for a control variable
$D$, we consider the following uncertain stochastic controlled system with
Markovian switching, for $t\in[0,T]$,
$$
\left\{\begin{array}{ll} dX(t)=&f(t, X(t), D(t),\zeta(t))dt+g(t,
X(t),D(t),\zeta(t))dB_t\\
&+h(t,X(t),D(t),\zeta(t))dC_t, \quad
X(0)=x,
\end{array}
\right. \eqno(2.1)
$$
where
$$
\begin{array}{ll}
&f=(f_1,\cdots,f_p)^\top: [0,T]\times{\mathbb R}^p\times
{\mathbb U}\times {\cal S}\rightarrow {\mathbb R}^p,\\
& g=(g_{kl})_{p\times m}: [0,T]\times{\mathbb R}^p\times{\mathbb U}\times {\cal S}\rightarrow {\mathbb
R}^{p\times m},\\
&h=(h_{kl})_{p\times n}: [0,T]\times{\mathbb R}^p\times{\mathbb U}\times {\cal S}\rightarrow {\mathbb
R}^{p\times n}.
\end{array}
$$

 Suppose that the cost functional of our control problem is as
follows
$$J(x, i; D)=E\left\{\int_0^T\Phi_1(t, X(t), D(t),\zeta(t))dt+\Phi_2(T,X(t), \zeta(T))\right\}. \eqno(2.2)$$
Set
$$
\begin{array}{ll}
{\cal U}\stackrel{\triangle}{=}&\left\{D: [0,T]\times\Gamma\times
\Omega\rightarrow {\mathbb U}\ |\ D\ {\rm is\ measurable\ and
\ for\ each }\right.\\
& \gamma\in\Gamma, D(\gamma)\ \mbox{is}\ \mathbb A \mbox{-adapted,\ and\ satisfies}\\
&\quad
E\left[\int_0^T\Phi_1^+(t,X(t), D(t),\zeta(t))dt\right]<+\infty,\\
&\quad
E\left[\int_0^T\Phi_2^+(T,X(t), \zeta(t))dt\right]<+\infty,\\
&\left.\quad E\left[\int_0^T|D(t)|^2dt\right]<+\infty
\right\}.
\end{array}
$$
Here,
$a^+\stackrel{\triangle}{=}\max(a,0),a^-\stackrel{\triangle}{=}\max(-a,0)$
for every $a\in[-\infty,\infty]$,  and $|\cdot|$ denotes the norm of
a vector or matrix.

   Under certain
assumptions (which will be specified below), for any $D(\cdot)\in
{\cal U}$, Eq. (2.1) admits a unique solution, and
the cost (2.2) is well-defined. The optimal control problem can be
stated as follows.

\vskip12pt

{\noindent \bf Problem 2.8.} Select an admissible control
$\hat{D}\in {\cal U}$ that minimizes $J(x,y,i;D))$ and find a value
function $V$ defined by $V(x, i)\stackrel{\triangle}{=}\inf_{D\in
{\cal U}}J(x, i; D).$ The control $\hat{D}$ is called an optimal
control.

\vskip12pt

 Next, we need to consider a family of optimal control problems
with different initial times and states along a given state
trajectory in order to apply the optimal principle for uncertain stochastic optimal control. For
this, we set up the following framework. For any $(t,
x)\in[0,T]\times {\mathbb R}^m$, consider the state equation,
$s\in[t,T]$,
$$
\left\{\begin{array}{ll}
 dX(s)=&f(s, X(s), D(s),\zeta(s))ds+g(s,
X(s),D(s),\zeta(s))dB_s\\
&+h(s,
X(s),D(s),\zeta(s))dC_s, \quad
X(t)=x.
\end{array}
\right.
\eqno(2.3)
$$

The corresponding cost functional is

$$J(t,x, i; D)=E\left[\int_t^T\Phi_1(s, X(s), D(s),\zeta(s))ds+\Phi_2(T, X(T), \zeta(T))\right]. \eqno(2.4)$$
Define the value function $V(t, x, i)$
by $V(t, x, i)\stackrel{\triangle}{=}\inf_{D\in {\cal U}}J(t,
x,i;D).$

For our aim, we make the following hypothesis.

\vskip12pt

{\noindent \bf Assumption 2.9.}  (1) Assume that coefficients $f,
g$ and $h$ are locally Lipschitzian continuous, that is,  there exists one nonnegative function $H_N(D)$, where
$H_N(D)$ is possibly dependent of the control variable $D$, such that the following inequality holds for $\hbar=f,g, h$,
$$
|\hbar(t, x, D, i)-\hbar(t,y, D,i)|<H_N(D)|x-y|,
$$
 for all $(t, D,
i)\in [0,T]\times {\cal U}\times {\cal S}$ and those
$x,y\in{\mathbb R}^p$ with $|x|\vee |y|\leq N.$ Moreover, for $g, h$, the
special growth conditions hold, that is, there is  a constant
$L_N$ such that
$$|g(t,x,D,i)|\vee|h(t,x,D,i)|<L_N(1+|x|+\tilde{H}(x)|D|),$$
for all $(t, D, i)\in [0,T]\times {\cal U}\times {\cal S}$ and $|x|\leq N$, where
the nonnegative function $\tilde{H}$ is locally bounded.

 (2) For the cost  functions  of the control problem $\Phi_l,l=1,2$, there exists a constant
 $K>0$,
 for all $i\in{\cal S}$,
 $$
 \Phi_1(t,x, D,i)\geq- K(1+\phi_1(x)|D|^2),\quad \Phi_2(t,x,i)\geq
- K(1+\phi_2(x)),
 $$
 where the nonnegative functions $\phi_l,l=1,2$, are continuous in
 $x$.

\vskip12pt

 Similar to the discussion for fuzzy differential equations in Fei \cite{FeiF}, it is known that under
Assumption 2.9 (1), Eq. (2.1) has a unique
continuous solution $X(t)=X(t,D(t),\zeta(t))$ on $t\in[0,T]$ for each $D\in {\cal U}$.

\vskip12pt

{\noindent \bf Remark 2.10.} We will see that Assumption 2.9 holds
for the linear system and special cost functions $\Phi_1$ and
$\Phi_2$ in the next section. For more general systems, we will
study the existence and uniqueness of solutions to uncertain stochastic
controlled systems, and the properties of value functions of our
control problems in future.

\vskip12pt

In order to solve our problem, now we generalize It\^o-Liu lemma for an uncertain stochastic process with Markovian swithcing as
follows. For sake of simplicity, we drop the arguments whenever
convenient, and denote
$$f(t)=f(t,X(t),D(t),\zeta(t)), g(t)=g(t,X(t),D(t),\zeta(t)), h(t)=h(t,X(t),D(t),\zeta(t)).$$

\vskip12pt

{\noindent\bf Lemma 2.11.}  {\sl
 Suppose that the state $X(t)=(X_1(t), \cdots, X_p(t))^\top$ of an uncertain stochastic system
satisfies Eq. (2.1). For
$\forall\ i\in {\cal S}, k,l=1, \cdots, p$, let $\frac{\partial F}{\partial t}(t, x, i)$, $\frac{\partial F}{\partial x_k}(t, x, i)$ and $\frac{\partial^2 F}{\partial x_kx_l}(t, x, i)$ be continuously functions, where $x=(x_1, \cdots, x_p)^\top$. Then we have
$$
\begin{array}{ll}
&{d}F(t, X(t), \zeta(t))\\
&=\frac{\partial F}{\partial t}(t, X(t), \zeta(t)){d }t+\sum\limits_{k=1}^p\frac{\partial F}{\partial x_k}(t, X(t),\zeta(t)){d}X_{k}(t)\\
&\quad+\frac{1}{2}\sum\limits_{k=1}^p\sum\limits_{l=1}^p\frac{\partial^2 F}{\partial x_k\partial x_l}(t, X(t),\zeta(t)){d}X_{k}(t){d}X_{l}(t)\\
&\quad+\sum\limits_{j=1}^Sq_{\zeta(t)j}F(t, X(t), j)dt+d M_t^F,
\end{array}
$$
where  the uncertain stochastic process $(M^F_t)_{0\leq t\leq T}$  is a real-valued
martingale relative to the filtered probability space $({\Omega, {\cal A}, {\mathbb A}, {\mathbb P}})$, ${d}B_{t}^k{d}B_{t}^l=\delta_{kl}{d}t, {d}B_{t}^k{d}t={d}C_{t}^\imath{d}C_{t}^\jmath={d}C_{t}^\imath{d}t={d}B_{t}^k{d}C_{t}^\imath=0,$ for
 $k,l=1,\cdots, m, \imath, \jmath=1, \cdots, n.$ Here, $\delta_{kl}$ is defined in Theorem 2.6.
 }

\vskip12pt

\Pr It is well known that almost every sample path of $\zeta(t)$ is
a right-continuous step function. It is useful to recall that a
continuous-time Markov chain $\zeta(t)$ can be represented as a
stochastic integral with respect to a Poisson random measure (cf.
Lemma 3 in Skorohod \cite{Sk}):
$$d\zeta(t)=\int_{\mathbb R}\hbar(X(t),\zeta(t),z)\nu(dz\times dt),$$
where $\nu(dz\times dt)$ is a Poisson random measure with intensity
$ m(dz)\times dt$ in which $m$ is a Lebesgue measure on $\mathbb R$.

Let $0<\tau_1<\cdots<\tau_n<t$ be all the times when $\zeta(t)$ has
a jump. From Theorem 2.6, we get
$$
\begin{array}{ll}
&F(\tau_1, X(\tau_1), \ell)-F(0, X(0), \ell)\\
&=\int_0^{\tau_1}\frac{\partial F}{\partial t}(s, X(s),\zeta(s)){d }s+\int_0^{\tau_1}\sum\limits_{k=1}^p\frac{\partial F}{\partial x_k}(s, X(s), \zeta(s)){d}X_{k}(s)\\
&\quad+\frac{1}{2}\int_0^{\tau_1}\sum\limits_{k=1}^p\sum\limits_{l=1}^p\frac{\partial^2 F}{\partial x_k\partial x_l}(s, X(s), \zeta(s)){d}X_{k}(s){d}X_{l}(s),
\end{array}
$$
$$
\begin{array}{ll}
&F(\tau_{k+1}, X(\tau_{k+1}), \ell)-F(\tau_k, X(\tau_k), \ell)\\
&=\int_{\tau_k}^{\tau_{k+1}}\frac{\partial F}{\partial t}(s, X(s),\zeta(s)){d }s+\int_{\tau_{k}}^{\tau_{k+1}}\sum\limits_{k=1}^p\frac{\partial F}{\partial x_k}(s, X(s), \zeta(s)){d}X_{k}(s)\\
&\quad+\frac{1}{2}\int_{\tau_{k}}^{\tau_{k+1}}\sum\limits_{k=1}^p\sum\limits_{l=1}^p\frac{\partial^2 F}{\partial x_k\partial x_l}(s, X(s),\zeta(s)){d}X_{k}(s){d}X_{l}(s),
\end{array}
$$
$$
\begin{array}{ll}
&F(t, X(t),  \ell)-F(\tau_n, X(\tau_n),  \ell)\\
&=\int_{\tau_n}^t\frac{\partial F}{\partial t}(s, X(s),\zeta(s)){d }s+\int_{\tau_n}^t\sum\limits_{k=1}^p\frac{\partial F}{\partial x_k}(s, X(s), \zeta(s)){d}X_{k}(s)\\
&\quad+\frac{1}{2}\int_{\tau_n}^t\sum\limits_{k=1}^p\sum\limits_{l=1}^p\frac{\partial^2 F}{\partial x_k\partial x_l}(s, X(s), \zeta(s)){d}X_{k}(s){d}X_{l}(s).
\end{array}
$$

Substituting $\ell=i$ in the first equation, $\ell=\zeta_{\tau_k}$ in the
second, and $\ell=\zeta_{\tau_n}$ in the third and adding them over
$k$ from 1 to $n+1$, we obtain
$$
 \begin{array}{ll}
 &F(t, X(t), \zeta(t))-F(0, X(0), i)\\
&=\int_0^t\frac{\partial F}{\partial t}(s, X(s),\zeta(s)){d }s+\int_0^t\sum\limits_{k=1}^p\frac{\partial F}{\partial x_k}(s, X(s),\zeta(s)){d}X_{k}(s)\\
&\quad+\frac{1}{2}\int_0^t\sum\limits_{k=1}^p\sum\limits_{l=1}^p\frac{\partial^2 F}{\partial x_k\partial x_l}(s, X(s),  \zeta(s)){d}X_{k}(s){d}X_{l}(s)\\
&\quad+\sum\limits_{k=1}^n[F(\tau_{k}, X(\tau_{k}), \zeta(\tau_k))-F(\tau_{k}, X(\tau_{k}),
\zeta(\tau_k-))]\\
& =\int_0^t\frac{\partial F}{\partial t}(s, X(s),\zeta(s)){d }s+\int_0^t\sum\limits_{k=1}^p\frac{\partial F}{\partial x_k}(s, X(s), \zeta(s)){d}X_{k}(s)\\
&\quad+\frac{1}{2}\int_0^t\sum\limits_{k=1}^p\sum\limits_{l=1}^p\frac{\partial^2 F}{\partial x_k\partial x_l}(s, X(s), \zeta(s)){d}X_{k}(s){d}X_{l}(s)\\
&\quad +\int_0^t\int_{\mathbb R}[F(s, X(s),  \zeta(s)+\hbar(X(s),\zeta(s),z))-F(s, X(s),\zeta(s))]m(dz)ds\\
&\quad +\int_0^t\int_{\mathbb R}[F(s, X(s), \zeta(s)+\hbar(X(s),\zeta(s),z))-F(s, X(s), \zeta(s))]\mu(dz\times ds),
\end{array}
$$
where $\mu(dz\times dt)=\nu(dz\times dt)-m(dz)\times dt$ is a
martingale measure relative to the filtered probability space $(\Omega, {\cal A}, {\mathbb A}, {\mathbb P})$ . Similar to the discussion of Lemma 3 (p. 104) in
Skorohod \cite{Sk}, noting that
$$
 \begin{array}{ll}
 &\int_0^t\int_{\mathbb R}[F(s, X(s),\zeta(s)+\hbar(X(s),\zeta(s),z))-F(s, X(s), \zeta(s))]m(dz)ds\\
&\quad=\sum\limits_{j=1}^Sq_{\zeta(s)j}F(s,X(s), j),
\end{array}
$$
we get the claim by differentiating the above equation. Thus, the
proof is complete. $\Box$

\vskip12pt

In Yong and Zhou \cite{YZ}, they discussed the principle of optimality and the HJB equation for stochastic optimal control. Next,
we similarly explore those for uncertain stochastic optimal control. First, we derive the principle of optimality for uncertain stochastic optimal control (2.4) with the constraint (2.3).

\vskip12pt

{\noindent\bf Theorem 2.12. } (Principle of optimality) {\sl Let assumption 2.9 hold. For any $(t, x)\in [0,T]\times {\mathbb R}^p$, we have

$$V(t, x,i)=\inf\limits_{D\in {\cal U}}E\left[\int_t^{\hat t}\Phi_1(s, X(s), D(s), \zeta(s))ds+V(\hat t,X(\hat t), \zeta(\hat t))\right],\ 0\leq t\leq \hat t\leq T. \eqno(2.5)$$

}

\vskip12pt

\Pr  From Assumption 2.9 (2), for $u\in{\cal U}$, $0\leq t\leq \hat t\leq T$, we get
$$
\begin{array}{ll}
&0\geq -E\left[\int_t^{\hat t}\Phi^-_1(s,X(s),D(s), \zeta(s))ds\right]\\
&\geq -
E\left[\int_t^{\hat t}(1+\phi_1(X(s))|D(s)|^2)ds\right]>-\infty,
\end{array}
$$
and
$$0\leq E\left[\int_t^{\hat t}\Phi^+_1(s,X(s),D(s),\zeta(s))ds\right]<+\infty,
$$
which shows
$$\left| E\left[\int_t^{\hat t}\Phi_1(s,X(s),D(s),\zeta(s))ds\right]\right|<+\infty.
$$
Likewise, we can obtain
$$|
E[\Phi_2(T,X(t),\zeta(T))]|<+\infty.
$$
Thus, we show the right hand of the identity (2.5) is well defined.

Next, we denote the right side of (2.5) by $\widetilde{V}(t, x, i)$. It follows from the definition
of $V(t, x,i)$ that, for any $D\in{\cal U}$,
$$
\begin{array}{ll}
V(t, x,i)\leq &E\left[\int_t^{\hat t}\Phi_1(s, X(s), D(s), \zeta(s))ds\right.\\
&\quad\left.+\int_{\hat t}^T\Phi_1(s, X(s), D(s), \zeta(s))ds+\Phi_2(T, X(T),\zeta(T))\right].
\end{array}
$$
Since the uncertain processes
$dC_s (s\in [t, \hat t))$ and $dC_s (s\in [\hat t, T ])$ are independent, we know that
 $$\int_t^{\hat t}\Phi_1(s, X(s), D(s), \zeta(s))ds\ \mbox{and}\ \int_{\hat t}^T\Phi_1(s, X(s), D(s), \zeta(s))ds$$
 are independent relative to the uncertain space $(\Gamma, {\cal L}, {\cal M})$.  Thus, by Theorem 1 in Zhu \cite{Zhu}, we get
 $$
\begin{array}{ll}
V(t, x,i)\leq &E\left\{\int_t^{\hat t}\Phi_1(s, X(s), D(s), \zeta(s))ds\right.\\
&\quad\left.+E\left[\int_{\hat t}^T\Phi_1(s, X(s), D(s), \zeta(s))ds+\Phi_2(T, X(T),\zeta(T))\right]\right\}.
\end{array}
$$
Taking the infimum for above inequality with respect to $D\in {\cal U}$, we have $V(t, x, i)\leq \widetilde{V}(t,x,i).$ On the other hand, for $\forall D\in {\cal U}$, we derive
$$
\begin{array}{ll}
&E\left[\int_t^T\Phi_1(s, X(s), D(s), \zeta(s))ds+\Phi_2(T,X(T), \zeta(T))\right]\\
&=E\left\{\int_t^{\hat t}\Phi_1(s, X(s), D(s), \zeta(s))ds\right.\\
&\quad\left.+E\left[\int_{\hat t}^T\Phi_1(s, X(s), D(s), \zeta(s))ds+\Phi_2(T, X(T),\zeta(T))\right]\right\}\\
&\geq E\left[\int_t^{\hat t}\Phi_1(s, X(s), D(s), \zeta(s))ds+V(\hat t, X(\hat t), \zeta(\hat t))\right]\geq \widetilde{V}(t, x, i)
\end{array}
$$
which shows $V(t, x, i)\geq \widetilde{V}(t,x,i).$  Consequently, $V(t, x, i)= \widetilde{V}(t,x,i).$  Thus, the proof of the theorem is complete. $\Box$

%%%%%%%%%%%%%%%%%%%%%%%%%

Let $C^{1,2}([0,T]\times{\mathbb R}^p; {\mathbb R })$ denote
 all functions $V(t,x,i)$ on $[0,T]\times{\mathbb R}^p$ which are continuously differentiable in $t$, continuously
twice differentiable in $x$ for each $i\in{\cal S}$. If
${V}(\cdot,i)\in C^{1,2}([0,T]\times{\mathbb R}^p;
{\mathbb R })$, define operators $L_i(D){V}$, for each $i\in{\cal
S}$, by
$$
 \begin{array}{ll}
&L_i(D){V}(t,x, i) \\
&\stackrel{\triangle}{=}\frac{\partial V}{\partial t}(t, x, i)+ \sum\limits_{k=1}^p\frac{\partial V}{\partial x_k}(t,x,i)f_{k}(t,x,D,i)\\
&\quad+\frac{1}{2}\sum\limits_{k=1}^p\sum\limits_{l=1}^p\frac{\partial^2 V}{x_{k}x_l}(t,x,i)\sum\limits_{\imath=1}^m g_{k\imath }(t,x,D,i)g_{l\imath }(t,x,D,i).
\end{array}
$$

In what follows, we give the optimal equation of the optimal
control problem.

\vskip12pt

{\noindent \bf Theorem 2.13.} {\sl Let Assumption 2.9 hold. Then $V$ is a solution of the following terminal problem of a Hamilton-Jacobi-Bellman (HJB) equation
$$\inf\limits_{D\in {\cal U}}\{L_i(D)V(t,x, i)+\Phi_1(t,x,D,i)\}
=\lambda_i V(t,x,i)-\sum\limits_{j\in{\cal
S}\backslash\{i\}}q_{ij}V(t,x,j)\eqno(2.6)$$
 with the terminal
condition $V(T, x,i)=\Phi_2(T,x,i).$
 }
\vskip12pt

\Pr  Let $a$ and $b$ satisfy $0<a<|X(t)|<b<+\infty$. We
define the first hitting times $\tau_a\stackrel{\triangle}{=}\inf\{t\geq 0:
|X(t)|=a\}, \tau_b\stackrel{\triangle}{=}\inf\{t\geq 0:
|X(t)|=b\}$, and
$\tau\stackrel{\triangle}{=}\tau_a\wedge\tau_b$.

For $\forall t, \hat t\in [0,T]$, from Lemma 2.11 we obtain
$$
\begin{array}{ll}
&V (\hat t\wedge\tau, X(\hat t\wedge\tau), \zeta(\hat t\wedge\tau))-V(t,x,i)\\
&=\int_t^{\hat t\wedge\tau}[L_{\zeta(s)}V(s,X(s),\zeta(s))+\sum\limits_{j=1}^Sq_{\zeta(s)j}V(s, X(s), j)]ds\\
&\quad+\int_t^{\hat t\wedge\tau}\sum\limits_{k=1}^p\sum\limits_{l=1}^mg_{kl}(s,X(s),D(s), \zeta(s))\frac{\partial V}{\partial x_k}(s,X(s),\zeta(s))dB_s^l\\
&\quad+\int_t^{\hat t\wedge\tau}\sum\limits_{k=1}^p\sum\limits_{l=1}^nh_{kl}(s,X(s),D(s),\zeta(s))\frac{\partial V}{\partial x_k}(s,X(s),\zeta(s))dC_s^l+ (M_{\hat t\vee \tau}^V-M_t^V).
\end{array}
\eqno(2.7)
$$

Since the uncertain stochastic process $(M_t^V)_{0\leq t\leq T}$ is a real-valued martingale relative to the filtered probability space $(\Omega, {\cal A}, {\mathbb A}, {\mathbb P})$, we easily know $E[M_{\hat t\vee\tau}^V]=E[M_t^V]$.
 Note that for all $i\in{\cal S}$, the functions $V_x(t,x,i)$ be bounded when $|x|$ is
bounded. Therefore, by Assumption 2.9 (1) and
$E[\int_0^T|D(t)|^2dt]<+\infty$, we easily deduce
$$
 \begin{array}{ll}
E\left|\int_t^{\hat t\wedge\tau}\sum\limits_{k=1}^p\sum\limits_{l=1}^mg_{kl}(s,X(s),D(s), \zeta(s))\frac{\partial V}{\partial x_k}(s,X(s),\zeta(s))dB_s^l\right|^2<+\infty.
\end{array}
$$
Hence, we get
$$
 \begin{array}{ll}
E\left[\int_t^{\hat t\wedge\tau}\sum\limits_{k=1}^p\sum\limits_{l=1}^mg_{kl}(s,X(s),D(s), \zeta(s))\frac{\partial V}{\partial x_k}(s,X(s),\zeta(s))dB_s^l\right]=0.
\end{array}
\eqno(2.8)
$$
On the other hand, similar to the discussion of Theorem 12.17 in Liu \cite{Liu4th}, it is easy to check that
$$
\begin{array}{ll}
E\left[\int_t^{\hat t\wedge\tau}\sum\limits_{k=1}^p\sum\limits_{l=1}^nh_{kl}(s,X(s),D(s),\zeta(s))\frac{\partial V}{\partial x_k}(s,X(s),\zeta(s))dC_s^l\right]=o(\hat t-t).
\end{array}
\eqno(2.9)
$$
Taking the expectation with respect to (2.7), together with (2.8) and (2.9), we have
$$
\begin{array}{ll}
&E[V (\hat t, X(\hat t), \zeta(\hat t))]=V(t,x,i)\\
&\quad +E\left[\int_t^{\hat t\wedge\tau}(L_{\zeta(s)}V(s,X(s),\zeta(s))+\sum\limits_{j=1}^Sq_{\zeta(s)j}V(s, X(s), j))ds\right]+o(\hat t-t).\\
&\end{array}
\eqno(2.10)
$$
 By using the principle of optimality (Theorem 2.12) and (2.10), we get
 $$
\begin{array}{ll}
0 =&\inf\limits_{D\in\cal U}E\left[\int_t^{\hat t\wedge\tau}\Phi_1(s, X(s), D(s), \zeta(s))ds\right.\\
&\quad+\left.\int_t^{\hat t\wedge\tau}(L_{\zeta(s)}V(s,X(s),\zeta(s))+\sum\limits_{j=1}^Sq_{\zeta(s)j}V(s, X(s), j))ds+o(\hat t-t)\right].
\end{array}
$$
Now letting $a\downarrow 0$ and $b\uparrow\infty$, we obtain
$\tau\rightarrow T.$ Therefore, applying the Monotone Convergence
Theorem, we have
$$
\begin{array}{ll}
0 =&\inf\limits_{D\in\cal U}E\left[\int_t^{\hat t}\Phi_1(s, X(s), D(s), \zeta(s))ds\right.\\
&\quad+\left.\int_t^{\hat t}(L_{\zeta(s)}V(s,X(s),\zeta(s))+\sum\limits_{j=1}^Sq_{\zeta(s)j}V(s, X(s), j))ds+o(\hat t-t)\right].
\end{array}
\eqno(2.11)
$$
Dividing Eq. (2.11) by $\hat t-t$, and letting $\hat t\rightarrow t$, we get Eq. (2.6). The terminal condition of Eq. (2.6) holds obviously. Thus, we complete the proof. $\Box$

 \vskip12pt

\section{ Uncertain random financial market set-up} \setzero

\vskip12pt
Applying the above results, we introduce a financial market with regime switching in this section.
 Let a filtered probability space  $(\Omega, {\cal A}, \mathbb{A}, {\mathbb P})$ hosting $m$-dimensional Wiener
process $(B_t)_{0\leq t\leq T}$ and an uncertain space $({\Gamma}, {\cal L}, {\cal M})$ hosting $n$-demensional canonical process $(C_t)_{0\leq t\leq T}$, respectively. The Wiener process
$(B_t)_{0\leq t\leq T}$ and the canonical process $(C_t)_{0\leq t\leq T}$ model risky uncertainties and behavior uncertainties in financial assets, respectively.
The regime switching is a Markov chain which has a strongly irreducible generator
$Q$ as in Section 2.

Assume that uncertain stochastic processes $P_0(t)$ and
$P_k(t),k=1,\dots,m,$ on $[0,T]$ represent the prices of the
riskless asset and the $m$ risky assets, respectively. They satisfy
the Markovian switching uncertain stochastic differential equations
 $$
 \begin{array}{ll}
dP_0(t)=&r_{\zeta(t)}(t)P_0(t)dt,\\
dP_k(t)=&\alpha^k_{\zeta(t)}(t)P_k(t)dt+\sigma_{\zeta(t)}^k(t)P_k(t)dB_t+\eta_{\zeta(t)}^k(t)P_k(t)dC_t,
\end{array}
$$
with initial prices $P_0(0)=1$ and $P^k(0)=p_k>0$, and initial state
$\zeta(0)=\zeta_0$. Here, $r_{\zeta(t)}(t)$,
$\alpha_{\zeta(t)}(t)=(\alpha_{\zeta(t)}^{1}(t),,\cdots,\alpha_{\zeta(t)}^{m}(t))^\top,$
$\sigma_{\zeta(t)}^{k}(t)=(\sigma_{\zeta(t)}^{k1}(t), \cdots,\sigma_{\zeta(t)}^{km}(t))$, and $\eta_{\zeta(t)}^{k}(t)=(\eta_{\zeta(t)}^{k1}(t), $ $\cdots, \eta_{\zeta(t)}^{k n}(t))\ (k=1,\cdots,m)$ are deterministic
and bounded interest rate, expected returns, random volatility
functions and uncertain volatility functions, respectively.
Denote
$$
\begin{array}{ll}
\sigma_{\zeta(t)}(t)=(\sigma_{\zeta(t)}^{1}(t),,\cdots,\sigma_{\zeta(t)}^{m}(t))^\top,
\eta_{\zeta(t)}(t)=(\eta_{\zeta(t)}^{1}(t),,\cdots,\eta_{\zeta(t)}^{m}(t))^\top.
\end{array}
$$

 We suppose that
$\Lambda_i(t)\stackrel{\triangle}{=}\sigma_i(t)\sigma_i(t)^\top,
i\in{\cal S}$, are positive definite, and the coefficients of the
market (i.e., $r,\alpha,\sigma, \eta$) depend on the regime of an
economy. We easily know $\Lambda_i(t)$ is deterministic and bounded.
The agent chooses a portfolio
$\pi=\{\pi(t)=(\pi^1(t),\cdots,\pi^m(t))^\top, t\in[0,T]\}$,
representing the fraction of wealth invested in each risky asset. We
need a technical condition to be satisfied. A portfolio vector
process is an uncertain stochastic vector process
$\pi$ such that
$E\left[\int_0^t\pi^\top(s)\Lambda_{\zeta(t)}\pi(s)ds\right]<+\infty$
for all $t\in[0,T]$. The
fraction of wealth invested in the riskless asset at time
$t\in[0,T]$ is then $1-\sum_{k=1}^m\pi^k(t)$. The investor also
chooses a consumption rate process $c=\{c(t),t\in[0,T]\}$: a
nonnegative uncertain stochastic process such that
$E\left[\int_0^tc^2(s)ds\right]<\infty$ for all $t\in[0,T]$.

Let the market price of risk of the market be
$\theta_{\zeta(t)}(t)\stackrel{\triangle}{=}\sigma^{-1}_{\zeta(t)}(t)(\alpha_{\zeta(t)}(t)-r_{\zeta}(t){\bf
1})$, where ${\bf 1}=(1,\cdots,1)^\top$. Now the  evolution of the
nominal wealth $W(t)$ at time $t$ can be written as
$$
\begin{array}{ll}
dW(t)=&r_{\zeta(t)}(t)W(t)dt+W(t)\pi^\top(t)\sigma_{\zeta(t)}(t)d[B_t+\theta_{\zeta(t)}(t)]\\
&\quad +W(t)\pi^\top(t)\eta_{\zeta(t)}(t)dC_t-c(t)dt.
\end{array}
\eqno(3.1)
$$

We introduce the utility functions $U_1(\cdot)$ and $U_2(\cdot)$ of consumption
and wealth, respectively, which are assumed to be twice
differentiable, strictly increasing, and concave. Moreover,
$U^\prime_l(0)=\infty,\ U^\prime_l(\infty)=0, l=1,2$. It is easily
to see that there exists a constant $K$ such that $U_l(y)\leq
K(1+y^2), l=1,2.$

Given $t\in [0, T]$, we now define the objective function which depends on the the market
regimes as follows
$$\mathbb J(t,x,i; \pi,c)=E\left[\int_t^Te^{-\beta (s-t)}U_1\left(c(s)\right)ds+e^{-\beta(T-t)}U_2\left(W(T)\right)\right],$$
 where for $s\in[t,T],$ the wealth $W(s)$ with $W(t)=x$ follows (3.1), $\beta$ is the utility discount rate, which may be
different from the risk-free rate. The control process
$D\stackrel{\triangle}{=}(\pi,c) $, that satisfies
$$E\left[\int_0^Te^{-\beta t}U_1^-\left(c(t)\right)dt+e^{-\beta T}U^-_2\left(W(T)\right)\right]<+\infty,$$
will be called an admissible control process. The set of all
admissible controls will be denoted by $\Xi$.  We define the
value function by
$$\mathbb V(t,x,i)=\sup\limits_{(\pi,c)\in{\Xi}}\mathbb J(t,x,i;\pi,c), \eqno(3.2)$$
which shows that the agent selects consumption and investment
processes in order to maximize the sum of his expected discounted
utilities from consumption and terminal wealth.

In what follows, we deduce the Hamilton-Jacobi-Bellman equation
satisfied by  optimal consumption and portfolio policies. For this,
denote operators $\mathbb L_i(D)\mathbb V$, for each $i\in{\cal
S}$, by
$$
 \begin{array}{ll}
&\mathbb L_i(D)\mathbb V(t,x, i) =\mathbb L(\pi,c)\mathbb V(t,x, i)\\
&\stackrel{\triangle}{=}\frac{1}{2}x^2\pi^\top\Lambda_i(t)\pi\mathbb V_{xx}+x\pi^\top\sigma_i(t)\theta_i(t)\mathbb V_x(t,x,i)\\
&\quad+r_i(t)x\mathbb V_x-c \mathbb V_{x}-\beta\mathbb V+\mathbb V_{t}(t,x,i).
\end{array}
$$
 We now give the equation of optimality for
the optimal policy as the main theorem in this section.

\vskip12pt

{\noindent \bf Theorem 3.1.} (Equation of optimality) {\sl   For each $i\in{\cal S},$ $\mathbb V(\cdot, i)$ is a solution of the following terminal problem of the HJB equation
$$\sup\limits_{(\pi,c)\in {\Xi}}\left[\mathbb{L}_i(\pi,c)\mathbb V(t,x, i)+U_1(c)\right]
=\lambda_i \mathbb V(t,x,i)-\sum\limits_{j\in{\mathbb{S}}\backslash\{i\}}q_{ij}\mathbb (t,x,j)\mathbb V(t,x, i)\eqno(3.3)$$
 with the terminal
condition $\mathbb V(T, x,i)=U_2(x).$
}
\vskip12pt

\Pr In terms of the state equation (3.1), we assume, in
Theorem 2.13,
$$
\begin{array}{ll}
&\Phi_1(t,x,\pi,c,i)=-e^{-\beta t}U_1(c),\\
&\Phi_2(T,x,i)=-e^{-\beta T}U_2(x),\\
&f(t,x,\pi,c,i)=r_i(t)x+x\pi^\top\sigma_i(t)\theta_i(t)-c,\\
&g(t,x,\pi,c,i)=x\pi^\top\sigma_i(t),\\
&h(t,x,\pi,c,i)=x\pi^\top\eta_i(t),
 \end{array}
$$
which show that the value
function $\mathbb V(t,x,i)$ in (3.2) equals to $-e^{\beta
t}V(t,x,i)$, where $V$ solves the equation (2.6). By the
definitions of utility functions $U_1(\cdot)$ and $U_2(\cdot)$, we
can take $\phi_1(x)=1$ and $\phi_2(x)=x$ in
Assumption 2.2 (2). Hence, the conditions of Theorem 2.13 are
fulfilled, moreover from the HJB equation (2.6), we derive that $\mathbb V(t,x,i)$
satisfies the Eq. (3.3). Thus, the proof is complete. $\Box$

\vskip12pt

{\noindent \bf Corollary 3.2.} {\sl   The optimal consumption and portfolio problem (3.2) has the optimal policy as follows
$$
\begin{array}{ll}
\hat{c}(t,x,\zeta(t)) = &\Psi_1(\mathbb V_x(t,x,\zeta(t))),\\
\hat{\pi}(t,x,\zeta(t))=&-\frac{\mathbb V_x(t,x,\zeta(t))}{x\mathbb V_{xx}(t,x,\zeta(t))}\left(\sigma_{\zeta(t)}^\top(t)\right)^{-1}\theta_{\zeta(t)}(t)
\end{array}
\eqno(3.4)
$$
where $\Psi_1(\cdot)$ is the unique inverse function of $U^\prime_1(\cdot)$.

}
\vskip12pt

\Pr Since the consumption policy can be solved by the concave
maximization in (3.3), we get the optimal feedback consumption
policy. In order to get the
portfolio policy, we make a quadratic maximization problem in (3.3). $\Box$

\vskip12pt

\section{Optimal  strategies for CRRA utility} \setzero

\vskip12pt

In this section, we discuss a special case of an investor with CRRA utility.  Let the utility
functions of an investor be
$U_l(z)=\frac{z^{1-\kappa}}{1-\kappa},\ \kappa>0,\ \kappa\neq 1,\ l=1,2.$
Thus, we obtain $\Psi_l(y)=y^{-1/\kappa}$, where $\Psi_l(\cdot)$ are
the inverse functions of $U^\prime_l(\cdot), l=1,2$. From (3.3) of
Theorem 3.1 and optimal policy $(\pi,c)$ of (3.4),  we
easily deduce that the value function $\mathbb V$ satisfies the HJB equation
as follows
$$
\begin{array}{ll}
&\mathbb V_t(t,x,i)-\beta\mathbb V(t,x,i)+r_i(t)x\mathbb V_x(t,x,i)+\frac{\kappa}{1-\kappa}\mathbb V_x(t,x,i)^{\frac{\kappa-1}{\kappa}}\\
&\quad-\frac{|\theta_i(t)\mathbb V_x(t,x,i)|^2}{2\mathbb V_{xx}(t,x,i)}=\lambda_i \mathbb V(t,x,i)-\sum\limits_{j\in{\cal S}\backslash
\{i\}}q_{ij}\mathbb V(t,x,j)
\end{array}
\eqno(4.1)
$$
with the terminal condition
$\mathbb V(T,x,i)=U_2(x)=\frac{1}{1-\kappa}x^{1-\kappa}$.

We guess the form of solution to the equation (4.1) as follows
$$\mathbb V(t,x,i)=\frac{1}{1-\kappa}A_i^\kappa(t)x^{1-\kappa},$$
which shows , for $ i=1,\cdots, S$,
$$
\kappa A_i^{\kappa-1}(t)A_i^\prime(t)-\rho_i(t) A_i^\kappa(t)+
\kappa A_i^{\kappa-1}(t)+\sum\limits_{j\in{\cal
S}\backslash\{i\}}q_{ij}A_j^\kappa(t)=0, \ A_i(T)=1,\eqno(4.2)
$$
where
$$
\begin{array}{ll}
\rho_i(t)\stackrel{\triangle}{=}\beta+\lambda_i-(1-\kappa)r_i(t)-\frac{1-\kappa}{2\kappa}|\theta_i(t)|^2.
\end{array}
$$

By the existence
and uniqueness theorem (cf., Hirsch et al. \cite{HSD}, p. 144), the system of equations (4.2) has the unique solution $\{A_i(t),
i=1,\cdots,S\}.$  From (3.4), we get
$$
\begin{array}{ll}
\hat{c}(t,x,i)=&\frac{x}{A_i(t)},\\
\hat{\pi}(t,x,i)=&\frac{1}{\kappa}\left(\sigma^\top_i(t)\right)^{-1}\theta_i(t).
\end{array}
\eqno(4.3)
$$

The formula (4.3) shows that an investor's consumption rate is
higher as he is wealthier. Also, the investor invests more in the
risky assets, the higher the market price of risk and the lower the
risks of the risky assets.

 From now on, for sake of convenience,  we only discuss the
following first-order nonlinear differential dynamic system (i.e.,
$S=2$)
$$
\begin{array}{ll}
&\kappa A_1^{\kappa-1}A_1^\prime-\rho_1 A_1^\kappa+ \kappa
A_1^{\kappa-1}+\lambda_1A_2^\kappa=0, \ A_1(T)=1,\\
& \kappa A_2^{\kappa-1}A_2^\prime-\rho_2 A_2^\kappa+ \kappa
A_2^{\kappa-1}+\lambda_2A_1^\kappa=0, \ A_2(T)=1.
\end{array}
\eqno(4.4)
$$

By the numerical solution to the equations (4.4), we give the
economic analysis of the related results.

Now let us consider the following numerical example. There are two
financial assets: a risk-free asset and a risky asset. The market coefficients
follow: $\alpha_1=0.15, \alpha_2=0.25$, $\sigma_1=0.25,
\sigma_2=0.6$, $r_1=0.05$, $r_2=0.01$, $\lambda_1=1.2,
\lambda_2=2.5, T=1$.

First, we consider a high risk-aversion investor with $\kappa=10$
and $\beta=0.07$. Using our numerical parameter values in (4.4),
from (4.3) we see how the optimal nominal consumption to wealth
ratio depends on time $t$. The results are shown in Figure 1.
Similarly, we consider a high risk-tolerance investor with
$\kappa=0.7$ and $\beta=0.8$. Using our numerical parameter values
in (4.4), from (4.3) we see how the optimal nominal consumption to
wealth ratio depends on time $t$. The results are shown in Figure 2.
From Figures 1-2, we see that an investor with high risk-aversion
increases his the optimal nominal consumption rate as the time gets
long, while his optimal consumption rate is larger during the state
of a good economy than during the state of a bad economy;
conversely, an investor with high risk tolerance decreases his the
optimal nominal consumption rate as the time gets long, while his
optimal consumption rate is smaller during the state of a good
economy than during the state of a bad economy.

%\vskip12pt

\begin{center}
\includegraphics[width=10cm]{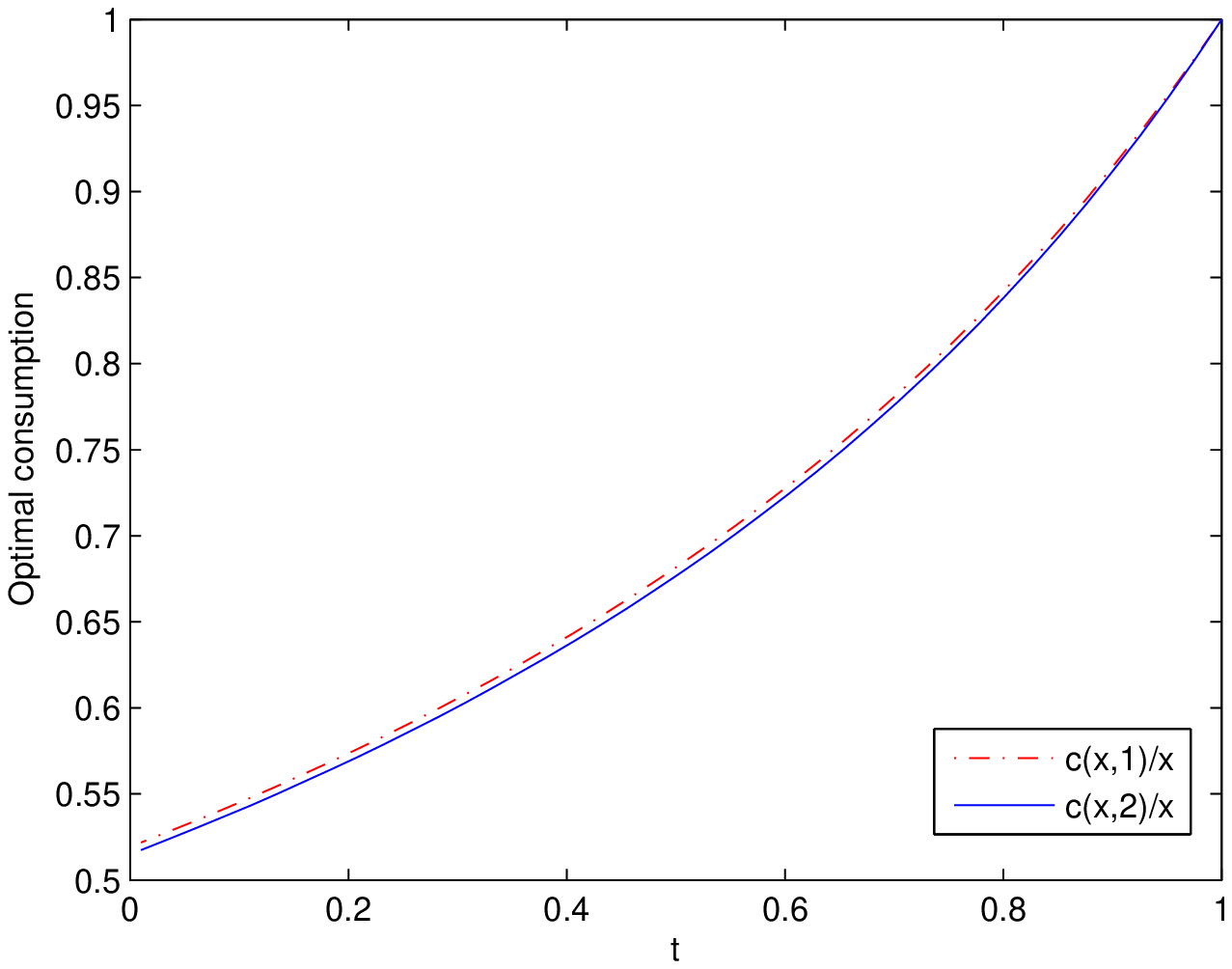}\\
{F{\scriptsize{IGURE}} 1.  An Investor with High Risk Aversion
($\kappa=10$ and $\beta=0.07$)}
\end{center}

 \vskip12pt

 \begin{center}
\includegraphics[width=10cm]{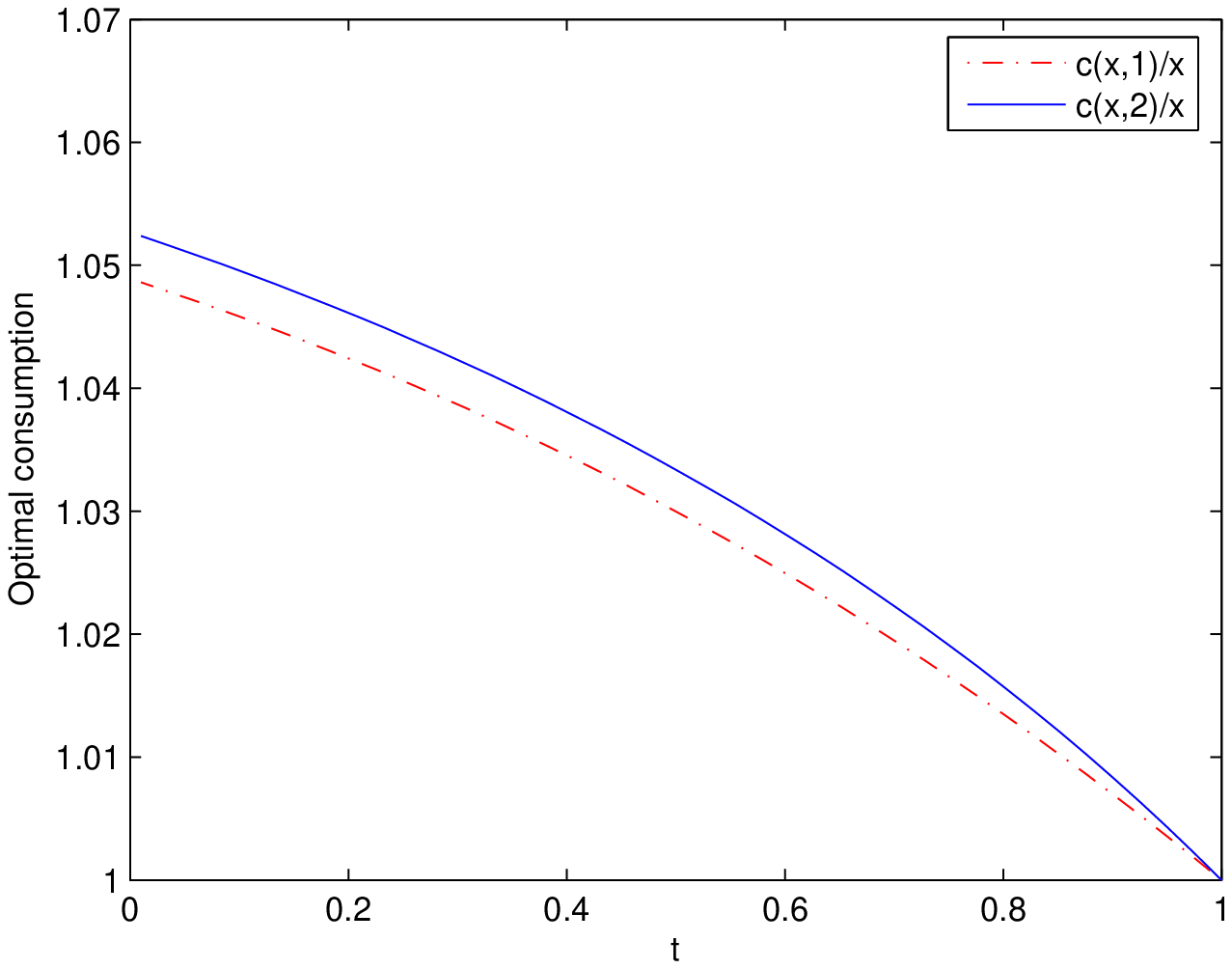}\\
{F{\scriptsize{IGURE}} 2.  An Investor with High Risk Tolerance
($\kappa=0.7$ and $\beta=0.8$)}
\end{center}

%\vskip12pt

\section{Conclusions and comments} \setzero

\vskip12pt

In this paper, we first describe an uncertain stochastic
control system with Markov-modulated. In order to build an optimality equation for
control laws, we generalize an It\^o-Liu formula for uncertain stochastic processes with Markovian switching.
 From this, the HJB equation of the
value function for the uncertain stochastic Markov-modulated control system is provided,
from which the optimal control law can be deduced. We think this is one of
the main contributions of this paper.

Following it, by the application of optimality equation for
optimal control law, we study the optimal consumption and
portfolio problem under an uncertain stochastic environment with regime switching, in which the prices  of the risky assets are driven by
risky uncertainty sources and Liu's uncertainty sources. An investor's utility comes from both the consumption and the wealth.
Employing the results on a Markovian-modulated uncertain stochastic control, we derive
the HJB equation satisfied by the optimal control law (consumption
and portfolio), from which we get the optimal consumption and portfolio polices. This is another contribution of our paper.
As far as we know, Zhou and Yin \cite{ZY} obtained explicit
solutions to classical stochastic control problems with regime
switching. Sotomayor and Cadenillas \cite{SC} presented the versions
of the HJB equation for classical stochastic control with regime
switching in infinite horizon. However, the framework  of our model
is different from theirs.

Finally, for an investor with CRRA utility functions, by using Corollary 3.2, we give an explicit formula for optimal policies.
 Through the numerical
simulation, an economic analysis of optimal policies is made. And we
find that the qualitative behavior of the consumption to wealth
ratio depends not only on the regime but also on the level of risk
aversion.

Besides, a stochastic financial market with inflation described in Bensoussan et al. \cite{BKS} can be extended to an uncertain stochastic market with both regime switching and inflation, which will be further studied in future.

\vskip12pt

%{\noindent\large\bf Acknowledgements:} The author is grateful to Professor Zidong Wang for his valuable comments and suggestions.

\end{document}